\documentstyle[11pt]{article}
\ifx\shlhetal\undefinedcontrolsequence\let\shlhetal\relax\fi
\input amssym.def
\input amssym.tex

\newcount\skewfactor
\def\mathunderaccent#1#2 {\let\theaccent#1\skewfactor#2
\mathpalette\putaccentunder}
\def\putaccentunder#1#2{\oalign{$#1#2$\crcr\hidewidth
\vbox to.2ex{\hbox{$#1\skew\skewfactor\theaccent{}$}\vss}\hidewidth}}
\def\name{\mathunderaccent\tilde-3 }

\newcommand{\bbQs}{{\Bbb  Q}^*}
\newcommand{\bbQ}{{\Bbb  Q}}

\newcommand{\V}{{\bf V}}
\newcommand{\constr}{{\bf L}}

\newcommand{\rest}{{\restriction}}

\newcommand{\lalse}{{}^{\textstyle \alpha\!>} 2}

\newcommand{\lh}{\ell{\rm g}}

\newcommand{\rng}{{\rm rang}}

\newcommand{\forces}{\Vdash} 
\newcommand{\lesdot}{\mathrel{\mathord{<}\!\!\raise 0.8
pt\hbox{$\scriptstyle\circ$}}}  

\newcommand{\QED}{\hfill\vrule width 6pt height 6pt depth 0pt\vspace{0.1in}} 
\newcommand{\Proof}{\noindent {\sc Proof} \hspace{0.2in}} 

\newcommand{\cH}{{\cal H}}
\newcommand{\hchi}{\big(\cH(\chi),\mathord{\in},\mathord{<^*_{\chi}}\big)} 

\newcommand{\pre}[2]{{}^{\textstyle #1}{#2}}
\newcommand{\conc}{^\frown\!}

\newcommand{\tbT}{\name{\bf T}}
\newcommand{\tA}{\name{A}}
\newcommand{\tbA}{\name{\bf A}}
\newcommand{\tbG}{\name{G}_{\bbQ}}
\newcommand{\tbI}{\name{\cal I}}
\newcommand{\BB}{{\frak B}}

\newcommand{\gra}{{\cal G}_\gamma(T,\bbQ)}
\newcommand{\cD}{{\cal D}}
\newtheorem{theorem}{Theorem}[section] 
\newtheorem{claim}{Claim}[theorem]
 
\newtheorem{proposition}[theorem]{Proposition} 
\newtheorem{corollary}[theorem]{Corollary} 
\newtheorem{definition}[theorem]{Definition}

\newtheorem{hypothesis}[theorem]{Hypothesis} 
\newtheorem{notation}[theorem]{Notation}

\title{On full Souslin trees}
\author{{\bf Saharon Shelah}\thanks{\ \ We thank the NSF for partially
supporting this research under grant \#144-EF67. Publication No 624.}\\
Institute of Mathematics\\
The Hebrew University of Jerusalem\\
91904 Jerusalem, Israel\\
and\\
Department of Mathematics\\
Rutgers University\\
New Brunswick, NJ 08854, USA\\
and\\
Mathematics Department\\
University of Wisconsin -- Madison\\
Madison, WI 53706, USA}

\begin{document}

\setcounter{page}{0}
\setcounter{section}{-1}

\maketitle

\begin{abstract}
In the present note we answer a question of Kunen (15.13 in \cite{Mi91})
showing (in \ref{first}) that 
\begin{quotation}
\noindent it is consistent that there are full Souslin trees.
\end{quotation}
\end{abstract}
\eject

\section{Introduction}
In the present paper we answer a combinatorial question of Kunen listed in
Arnie Miller's Problem List. We force, e.g.~for the first strongly
inaccessible Mahlo cardinal $\lambda$, a full (see \ref{full}(2))
$\lambda$--Souslin tree and we remark that the existence of such trees follows
from $\V=\constr$ (if $\lambda$ is Mahlo strongly inaccessible). This answers
\cite[Problem 15.13]{Mi91}. 

Our notation is rather standard and compatible with those of classical
textbooks on Set Theory. However, in forcing considerations, we keep the older
tradition that  
\begin{center}
{\em 
a stronger condition is the larger one.
}
\end{center}
We will keep the following conventions concerning use of symbols.
\begin{notation}
{\em
\begin{enumerate}
\item $\lambda,\mu$ will denote cardinal numbers and $\alpha,\beta,\gamma, 
\delta,\xi,\zeta$ will be used to denote ordinals. 
\item Sequences (not necessarily finite) of ordinals are denoted by $\nu$,
$\eta$, $\rho$ (with possible indexes).
\item The length of a sequence $\eta$ is $\lh(\eta)$.
\item For a sequence $\eta$ and an ordinal $\alpha\leq\lh(\eta)$, $\eta\rest
\alpha$ is the restriction of the sequence $\eta$ to $\alpha$ (so $\lh(\eta
\rest\alpha)=\alpha$). If a sequence $\nu$ is a proper initial segment of a
sequence $\eta$ then we write $\nu\vartriangleleft\eta$ (and $\nu
\trianglelefteq\eta$ has the obvious meaning).  
\item A tilde indicates that we are dealing with a name for an object in
forcing extension (like $\name{x}$).
\end{enumerate}
}
\end{notation}

\section{Full $\lambda$-Souslin trees}
A subset $T$ of $\lalse$ is an $\alpha$--tree whenever ($\alpha$ is a limit
ordinal and) the following three conditions are satisfied: 
\begin{itemize}
\item $\langle\rangle\in T$, if $\nu\vartriangleleft\eta\in T$ then $\nu\in
T$, 
\item $\eta\in T$ implies $\eta\conc\langle 0\rangle,\;\eta\conc\langle
1\rangle\in T$,\quad and 
\item for every $\eta\in T$ and $\beta<\alpha$ such that $\lh(\eta)\leq\beta$
there is $\nu\in T$ such that $\eta\trianglelefteq\nu$ and $\lh(\eta)=\beta$.
\end{itemize}
A $\lambda$--Souslin tree is a $\lambda$--tree $T\subseteq\pre{\lambda{>}}{2}$
in which every antichain is of size less than $\lambda$. 

\begin{definition}
\label{full}
\begin{enumerate}
\item For a tree $T\subseteq\lalse$ and an ordinal $\beta\leq\alpha$ we let 
\[T_{[\beta]}\stackrel{\rm def}{=} T\cap\pre{\beta}{2}\quad\mbox{ and }\quad
T_{[<\beta]}\stackrel{\rm def}{=}  T\cap\pre{\beta{>}}{2}.\]
If $\delta\leq\alpha$ is limit then we define
\[{\lim}_\delta T_{[<\delta]}\stackrel{\rm def}{=}\{\eta\in\pre{\delta}{2}:
(\forall\beta<\delta)(\eta\rest\beta\in T)\}.\]
\item An $\alpha$--tree $T$ is {\em full} if for every limit ordinal $\delta<
\alpha$ the set $\lim_\delta(T_{[<\delta]})\setminus T_{[\delta]}$ has at most
one element. 
\item An $\alpha$--tree $T\subseteq\lalse$ has {\em true height $\alpha$} if
for every $\eta\in T$ there is $\nu\in\pre{\alpha}2$ such that 
\[\eta\vartriangleleft\nu\quad\mbox{ and }\quad (\forall\beta<\alpha)(\nu\rest
\beta\in T).\]
\end{enumerate}
\end{definition}
We will show that the existence of full $\lambda$--Souslin trees is consistent
assuming the cardinal $\lambda$ satisfies the following hypothesis. 

\begin{hypothesis}
\label{hypot}
{\em 
\begin{description}
\item[(a)] $\lambda$ is strongly inaccessible (Mahlo) cardinal,
\item[(b)] $S\subseteq\{\mu<\lambda:\mu$ is a strongly inaccessible cardinal
$\}$ is a stationary set, 
\item[(c)] $S_0\subseteq\lambda$ is a set of limit ordinals,
\item[(d)] for every cardinal $\mu\in S$, $\diamondsuit_{S_0\cap \mu}$ holds
true. 
\end{description}
Further in this section we will assume that $\lambda$, $S_0$ and $S$ are as
above and we may forget to repeat these assumptions. 
}
\end{hypothesis}
Let as recall that the diamond principle $\diamondsuit_{S_0\cap \mu}$
postulates the existence of a sequence $\bar{\nu}=\langle\nu_\delta:\delta\in 
S_0\cap\mu\rangle$ (called {\em a $\diamondsuit_{S_0\cap\mu}$--sequence}) such
that $\nu_\delta\in \pre{\delta}{2}$ (for $\delta\in S_0\cap\mu$) and 
\[(\forall\nu\in\pre{\mu}{2})[\mbox{ the set }\{\delta\in S_0\cap\mu:\nu
\restriction\delta =\nu_\delta\}\mbox{ is stationary in }\mu].\]

Now we introduce a forcing notion $\bbQ$ and its relative $\bbQs$ which will
be used in our proof.

\begin{definition}
\begin{enumerate}
\item {\bf A condition} in $\bbQ$ is a tree $T\subseteq \lalse$ of a true
hight $\alpha=\alpha(T)<\lambda$ (see \ref{full}(3); so $\alpha$ is a limit
ordinal) such that $\|\lim_\delta (T_{[<\delta]})\setminus T_{[\delta]}\| \le
1$ for every limit ordinal $\delta<\alpha$, 

{\bf the order} on $\bbQ$ is defined by\qquad $T_1\leq T_2$ if and only if 

$T_1=T_2\cap \pre{\alpha(T_1){>}}{2}$ (so it is the end--extension order).
\item For a condition $T\in\bbQ$ and a limit ordinal $\delta<\alpha(T)$, let
$\eta_\delta(T)$ be the unique member of $\lim_\delta(T_{[<\delta]})\setminus
T_{[\delta]}$ if there is one, otherwise $\eta_\delta(T)$ is not defined.
\item Let $T\in\bbQ$. A function $f:T\longrightarrow\lim_{\alpha(T)}(T)$ is
called {\em a witness for $T$} if $(\forall\eta\in T)(\eta\vartriangleleft
f(\eta))$.
\item {\bf A condition} in $\bbQs$ is a pair $(T,f)$ such that $T\in\bbQ$ and
$f:T\longrightarrow\lim_{\alpha(T)}(T)$ is a witness for $T$,

{\bf the order} on $\bbQs$ is defined by\qquad $(T_1,f_1)\leq (T_2,f_2)$ if
and only if 

$T_1\leq_{\bbQ} T_2$ and $(\forall\eta\in T_1)(f_1(\eta )\trianglelefteq
f_2(\eta))$. 
\end{enumerate}
\end{definition}

\begin{proposition}
\label{trivial}
\begin{enumerate}
\item If $(T_1,f_1)\in\bbQs$, $T_1\leq_{\bbQ} T_2$ and
\begin{description}
\item[($*$)] either $\eta_{\alpha(T_1)}(T_2)$ is not defined or it does not
belong to $\rng(f_1)$
\end{description}
then there is $f_2:T_2\longrightarrow\lim_{\alpha(T_2)}(T_2)$ such that $(T_1,
f_1)\leq (T_2,f_2)\in\bbQs$.
\item For every $T\in \bbQ$ there is a witness $f$ for $T$.
\end{enumerate}
\end{proposition}

\Proof Should be clear. \QED

\begin{proposition}
\label{easy}
\begin{enumerate}
\item The forcing notion $\bbQs$ is $(<\lambda)$--complete, in fact any
increasing chain of length $<\lambda$ has the least upper bound in $\bbQs$. 
\item The forcing notion $\bbQ$ is strategically $\gamma$-complete for each
$\gamma<\lambda$. 
\item Forcing with $\bbQ$ adds no new sequences of length $<\lambda$. Since
$\|\bbQ\| =\lambda$, forcing with $\bbQ$ preserves cardinal numbers,
cofinalities and cardinal arithmetic. 
\end{enumerate}
\end{proposition}

\Proof 1)\ \ \ It is straightforward: suppose that $\langle (T_\zeta,f_\zeta):
\zeta<\xi\rangle$ is an increasing sequence of elements of $\bbQs$. Clearly we
may assume that $\xi<\lambda$ is a limit ordinal and $\zeta_1<\zeta_2<\xi\quad
\Rightarrow\quad \alpha(T_{\zeta_1})<\alpha(T_{\zeta_2})$. Let $T_\xi=
\bigcup\limits_{\zeta<\xi}T_\zeta$ and $\alpha=\sup\limits_{\zeta<\xi}\alpha(
T_\zeta)$. Easily, the union is increasing and the $T_\xi$ is a full
$\alpha$--tree. For $\eta\in T_\xi$ let $\zeta_0(\eta)$ be the first $\zeta<
\xi$ such that $\eta\in T_\zeta$ and let $f_\xi(\eta)=\bigcup\{f_\zeta(\eta):
\zeta_0(\eta)\leq\zeta<\xi\}$. By the definition of the order on $\bbQs$ we
get that the sequence $\langle f_\zeta(\eta):\zeta_0(\eta)\leq\zeta<\xi
\rangle$ is $\vartriangleleft$--increasing and hence $f_\xi(\eta)\in
\lim_\alpha(T_\xi)$. Plainly, the function $f_\xi$ witnesses that $T_\xi$ has
a true height $\alpha$, and thus $(T_\xi,f_\xi)\in\bbQs$. It should be clear
that $(T_\xi,f_\xi)$ is the least upper bound of the sequence $\langle
(T_\zeta,f_\zeta): \zeta<\xi\rangle$.

\noindent 2)\ \ \ For our purpose it is enough to show that for each ordinal
$\gamma<\lambda$ and a condition $T\in\bbQ$ the second player has a winning
strategy in the following game $\gra$. (Also we can let Player I choose
$T_\xi$ for $\xi$ odd.)
\begin{quotation}
\noindent The game lasts $\gamma$ moves and during a play the players, called
I and II, choose successively open dense subsets $\cD_\xi$ of $\bbQ$ and
conditions $T_\xi\in\bbQ$. At stage $\xi<\gamma$ of the game:\\
Player I chooses an open dense subset $\cD_\xi$ of $\bbQ$ and\\
Player II answers playing a condition $T_\xi\in \bbQ$ such that
\[T\leq_{\bbQ} T_\xi,\quad (\forall\zeta<\xi)(T_\zeta\leq_{\bbQ} T_\xi),\quad
\mbox{ and}\quad T_\xi\in \cD_\xi.\]
The second player wins if he has always legal moves during the play. 
\end{quotation}
Let us describe the winning strategy for Player II. At each stage $\xi<\gamma$
of the game he plays a condition $T_\xi$ and writes down on a side a function
$f_\xi$ such that $(T_\xi,f_\xi)\in\bbQs$. Moreover, he keeps an extra
obligation that $(T_\zeta,f_\zeta)\leq_{\bbQs} (T_\xi,f_\xi)$ for each
$\zeta<\xi<\gamma$.\\
So arriving to a non-limit stage of the game he takes the condition
$(T_\zeta,f_\zeta)$ he constructed before (or just $(T,f)$, where $f$ is a
witness for $T$, if this is the first move; by \ref{trivial}(2) we can always
find a witness). Then he chooses $T^*_\zeta\geq_{\bbQ} T_\zeta$ such that
$\alpha(T^*_\zeta)=\alpha(T_\zeta)+\omega$ and $(T^*_\zeta)_{[\alpha(T_\zeta
)]}=\lim_{\alpha(T_\zeta)}(T_\zeta)$. Thus $\eta_{\alpha(T_\zeta)}(T^*_\zeta)$
is not defined. Now Player II takes $T_{\zeta+1}\geq_{\bbQ} T^*_\zeta$ from
the open dense set $\cD_{\zeta+1}$ played by his opponent at this stage. 
Clearly $\eta_{\alpha(T_\zeta)}(T_{\zeta+1})$ is not defined, so Player II may
use \ref{trivial}(1) to choose $f_{\zeta+1}$ such that $(T_\zeta,f_\zeta)
\leq_{\bbQs}(T_{\zeta+1},f_{\zeta+1})\in\bbQs$.\\
At a limit stage $\xi$ of the game, the second player may take the least upper
bound $(T^\prime_\xi,f^\prime_\xi)\in\bbQs$ of the sequence $\langle (T_\zeta,
f_\zeta): \zeta<\xi\rangle$ (exists by 1))  and then apply the procedure
described above. 

\noindent 3)\ \ \ Follows from 2) above. \QED

\begin{definition}
\label{generic}
Let $\tbT$ be the canonical $\bbQ$--name for a generic tree added by forcing
with $\bbQ$:
\[\forces_{\bbQ}\tbT=\bigcup\{T:T\in\tbG\}.\]
\end{definition}
It should be clear that $\tbT$ is (forced to be) a full $\lambda$--tree. The
main point is to show that it is $\lambda$--Souslin and this is done in the
following theorem. 

\begin{theorem}
\label{first}
$\forces_{\bbQ}$`` $\tbT$ is a $\lambda$--Souslin tree''.
\end{theorem}

\Proof Suppose that $\tA$ is a $\bbQ$--name such that
\[\forces_{\bbQ}\mbox{`` }\tA\subseteq\tbT\mbox{ is an antichain '',}\]
and let $T_0$ be a condition in $\bbQ$. We will show that there are $\mu<
\lambda$ and a condition $T^*\in\bbQ$ stronger than $T_0$ such that
$T^*\forces_{\bbQ}$ ``$\tA\subseteq\tbT_{[<\mu]}$ '' (and thus it forces that
the size of $\tA$ is less than $\lambda$). 

Let $\tbA$ be a $\bbQ$--name such that
\[\forces_{\bbQ}\mbox{`` }\tbA=\{\eta\in\tbT: (\exists\nu\in\tA)(\nu
\trianglelefteq\eta)\mbox{ or }\neg(\exists\nu\in\tA)(\eta\trianglelefteq\nu)
\}\mbox{ ''.}\]
Clearly, $\forces_{\bbQ}$ ``$\tbA\subseteq\tbT$ is dense open''.\\
Let $\chi$ be a sufficiently large regular cardinal ($\beth_7(\lambda^+)^+$
is enough). 

\begin{claim}
\label{cl1}
There are $\mu\in S$ and $\BB\prec\hchi$ such that:
\begin{description}
\item[(a)] $\tA,\tbA,S, S_0,\bbQ,\bbQs,T_0\in\BB$,
\item[(b)] $\|\BB\| =\mu$\quad and\quad $\pre{\mu{>}}{\BB}\subseteq\BB$,
\item[(c)] $\BB\cap\lambda =\mu$.
\end{description}
\end{claim}

\noindent{\em Proof of the claim:}\ \ \ First construct inductively an
increasing continuous sequence $\langle\BB_\xi: \xi<\lambda\rangle$ of
elementary submodels of $\hchi$ such that $\tA,\tbA,S, S_0,\bbQ,\bbQs,T_0\in
\BB_0$ and for every $\xi<\lambda$
\[\|\BB_\xi\|=\mu_\xi<\lambda,\quad\BB_\xi\cap\lambda\in\lambda,\quad\mbox{
and }\quad\pre{\mu_\xi{\geq}}{\BB_\xi}\subseteq\BB_{\xi+1}.\]
Note that for a club $E$ of $\lambda$, for every $\mu\in S\cap E$ we have 
\[\|\BB_\mu\| =\mu,\quad\pre{\mu{>}}{\BB_\mu}\subseteq\BB_\mu,\quad\mbox{ and
}\quad \BB\cap\lambda =\mu.\]
\smallskip

Let $\mu\in S$ and $\BB\prec\hchi$ be given by \ref{cl1}. We know that
$\diamondsuit_{S_0\cap\mu}$ holds, so fix a $\diamondsuit_{S_0\cap
\mu}$--sequence $\bar{\nu}=\langle\nu_\delta :\delta\in S_0\cap\mu\rangle$. 

Let 
\[\begin{array}{ll}
\tbI\stackrel{\rm def}{=}\big\{T\in\bbQ: &T\mbox{ is incompatible (in $\bbQ$)
with }T_0\ \ \mbox{ or:}\\
\ &T\geq T_0\mbox{ and } T\mbox{ decides the value of }\tbA\cap\pre{\alpha(T)
  {>}}{2}\ \mbox{ and}\\
\ &(\forall\eta\in T)(\exists\rho\in T)(\eta\trianglelefteq\rho\ \&\
  T\forces_{\bbQ}\rho\in\tbA)\big\}.
  \end{array}\]

\begin{claim}
\label{cl2}
$\tbI$ is a dense subset of $\bbQ$.
\end{claim}

\noindent{\em Proof of the claim:}\ \ \ Should be clear (remember
\ref{easy}(2)). 
\medskip

Now we choose by induction on $\xi<\mu$ a continuous increasing sequence
$\langle(T_\xi, f_\xi):\xi<\mu\rangle\subseteq \bbQs\cap\BB$.

\noindent{\sc Step:}\quad $i=0$\\
$T_0$ is already chosen and it belongs to $\bbQ\cap\BB$. We take any $f_0$
such that $(T_0,f_0)\in\bbQs\cap\BB$ (exists by \ref{trivial}(2)).

\noindent{\sc Step:}\quad limit $\xi$\\
Since $\pre{\mu{>}}{\BB}\subseteq\BB$, the sequence $\langle(T_\zeta,f_\zeta):
\zeta<\xi\rangle$ is in $\BB$. By \ref{easy}(1) it has the least upper bound
$(T_\xi,f_\xi)$ (which belongs to $\BB$).

\noindent{\sc Step:}\quad $\xi=\zeta+1$\\
First we take (the unique) tree $T^*_\xi$ of true height $\alpha(T^*_\xi)=
\alpha(T_\zeta)+\omega$ such that $T^*_\xi\cap\pre{\alpha(T_\zeta){>}}{2}=
T_\zeta$ and:\\
if $\alpha(T_\zeta)\in S_0$ and $\nu_{\alpha(T_\zeta)}\notin\rng(f_\zeta)$
then $(T^*_\xi)_{[\alpha(T_\zeta)]}=\lim_{\alpha(T_\zeta)}(T_\zeta)\setminus
\{\nu_{\alpha(T_\zeta)}\}$,\\
otherwise $(T^*_\xi)_{[\alpha(T_\zeta)]}=\lim_{\alpha(T_\zeta)}(T_\zeta)$.\\
Let $T_\xi\in\bbQ\cap\tbI$ be strictly above $T^*_\xi$ (exists by
\ref{cl2}). Clearly we may choose such $T_\xi$ in $\BB$. Now we have to define
$f_\xi$.  We do it by \ref{trivial}, but additionally we require that
\[\mbox{if }\quad\eta\in T_\xi\quad\mbox{ then }\quad(\exists\rho\in T_\xi)(
\rho\vartriangleleft f_\xi(\eta)\ \& \ T\forces_{\bbQ}\mbox{`` }\rho\in\tbA
\mbox{ ''}).\]
Plainly the additional requirement causes no problems (remember the definition
of $\tbI$ and the choice of $T_\xi$) and the choice can be done in $\BB$.
\medskip

There are no difficulties in carrying out the induction. Finally we let 
\[T_\mu\stackrel{\rm def}{=}\bigcup_{\xi<\mu}T_\xi\quad\mbox{ and }\quad
f_\mu=\bigcup_{\xi<\mu}f_\xi.\]
By the choice of $\BB$ and $\mu$ we are sure that $T_\mu$ is a $\mu$--tree. It
follows from \ref{easy}(1) that $(T_\mu,f_\mu)\in\bbQs$, so in particular 
the tree $T_\mu$ has enough $\mu$  branches (and belongs to $\bbQ$). 

\begin{claim}
\label{cl3}
For every $\rho\in\lim_\mu(T_\mu)$ there is $\xi<\mu$ such that
\[(\exists\beta<\alpha(T_{\xi+1}))(T_{\xi+1}\forces_{\bbQ}\mbox{`` }\rho\rest
\beta\in\tbA\mbox{ ''}).\]
\end{claim}

\noindent{\em Proof of the claim:}\ \ \ Fix $\rho\in\lim_\mu(T_\mu )$ and let 
\[S^*_\nu\stackrel{\rm def}{=}\{\delta\in S_0\cap\mu:\alpha(T_\delta )=\delta\ 
\mbox{ and }\ \nu_\delta=\rho\rest\delta\}.\]
Plainly, the set $S^*_\nu$ is stationary in $\mu$ (remember the choice of
$\bar{\nu}$). By the definition of the $T_\xi$'s (and by $\rho\in\lim_\mu
(T_\mu)$) we conclude that for every $\delta\in S^*_\nu$
\begin{quotation}
\noindent if $\eta_\delta(T_{\delta+1})$ is defined\quad then $\rho\rest\delta
\neq\eta_\delta(T_\mu)=\eta_\delta(T_{\delta+1})$.
\end{quotation}
But $\rho\restriction\delta =\nu_\delta$ (as $\delta\in S^*_\nu$). So look at
the inductive definition: necessarily for some $\rho^*_\delta\in T_\delta$ we
have $\nu_\delta =f_\delta(\rho^*_\delta )$, i.e.~$\rho\rest\delta=f_\delta(
\rho^*_\delta)$. Now, $\rho^*_\delta\in T_\delta=\bigcup\limits_{\xi<\delta}
T_\xi$ and hence for some $\xi(\delta)<\delta$, we have $\rho^*_\delta\in
T_{\xi(\delta )}$.  By Fodor lemma we find $\xi^*<\mu$ such that the set 
\[S'_\nu\stackrel{\rm def}{=}\{\delta\in S^*_\nu:\xi(\delta )=\xi^*\}\]
is stationary in $\mu$. Consequently we find $\rho^*$ such that the set
\[S^+_\nu\stackrel{\rm def}{=}\{\delta\in S'_\nu:\rho^*=\rho^*_\delta\}\]
is stationary (in $\mu$). But the sequence $\langle f_\xi(\rho^*):\xi^*\leq
\xi<\mu\rangle$ is $\trianglelefteq$--increasing, and hence the sequence
$\rho$ is its limit. Now we easily conclude the claim using the inductive
definition of the $(T_\xi,f_\xi)$'s.
\medskip

It follows from the definition of $\tbA$ and \ref{cl3} that
\[T_\mu\forces_{\bbQ}\mbox{`` }\tA\subseteq T_\mu\mbox{ ''}\]
(remember that $\tA$ is a name for an antichain of $\tbT$), and hence
\[T_\mu\forces_{\bbQ}\mbox{`` }\|\tA\|<\lambda\mbox{ ''},\]
finishing the proof of the theorem. \QED

\begin{definition}
A $\lambda$--tree $T$ is {\em $S_0$--full}, where $S_0\subseteq\lambda$, if
for every limit $\delta <\lambda$
\begin{quotation}
\noindent if $\delta\in\lambda\setminus S_0$ then $T_{[\delta]}=
\lim_\delta(T)$,\\
if $\delta\in S_0$ then $\|T_{[\delta]}\setminus\lim_\delta(T)\|\leq 1$.
\end{quotation}
\end{definition}

\begin{corollary}
Assuming Hypothesis \ref{hypot}:
\begin{enumerate}
\item The forcing notion $\bbQ$ preserves cardinal numbers, cofinalities and
cardinal arithmetic.
\item $\forces_{\bbQ}$ `` $\tbT\subseteq\pre{\lambda{>}}{2}$ is a
$\lambda$--Souslin tree which is full and even $S_0$--full ''.

[So, in $\V^{\bbQ}$, in particular we have:

for every $\alpha<\beta<\mu$, for all $\eta\in T\cap \pre{\alpha}{2}$ there is
$\nu\in T\cap\pre{\beta}{2}$ such that $\eta\vartriangleleft\nu$, and for a
limit ordinal $\delta<\lambda$, $\lim_\delta(T_{[<\delta]})\setminus 
T_{[\delta]}$ is  either empty or has a unique element (and then $\delta\in
S_0$).]
\end{enumerate}
\end{corollary}

\Proof By \ref{easy} and \ref{first}. \QED
\medskip

Of course, we do not need to force.

\begin{definition}
\label{sqdia}
Let $S_0,S\subseteq\lambda$. A sequence $\langle(C_\alpha,\nu_\alpha ):\alpha
<\lambda\mbox{ limit}\rangle$ is called {\em a squared diamond sequence for
$(S,S_0)$} if for each limit ordinal $\alpha<\lambda$
\begin{description}
\item[(i)]\ \ \ $C_\alpha$ a club of $\alpha$ disjoint to $S$,
\item[(ii)]\ \  $\nu_\alpha\in\pre{\alpha}{2}$,
\item[(iii)]\   if $\beta\in {\rm acc}(C_\alpha)$ then $C_\beta=C_\alpha\cap
\beta$ and $\nu_\beta\vartriangleleft \nu_\alpha$,
\item[(iv)]\ \  if $\mu\in S$ then $\langle \nu_\alpha:\alpha\in C_\mu\cap
S_0\rangle$ is a diamond sequence. 
\end{description}
\end{definition}

\begin{proposition}
Assume (in addition to \ref{hypot})
\begin{description}
\item[(e)] there exist a squared diamond sequence for $(S,S_0)$.
\end{description}
Then there is a $\lambda$--Souslin tree $T\subseteq \pre{\lambda{>}}2$ which
is $S_0$--full.
\end{proposition}

\Proof Look carefully at the proof of \ref{first}. \QED

\begin{corollary}
Assume that $\V=\constr$ and $\lambda$ is Mahlo strongly inaccessible. Then
there is a full $\lambda$--Souslin tree.
\end{corollary}

\Proof Let $S\subseteq\{\mu<\lambda: \mu\mbox{ is strongly inaccessible}\}$ be
a stationary non-reflecting set. By Beller and Litman \cite{BeLi80}, there is
a square $\langle C_\delta :\delta <\lambda\;\mbox{ limit}\rangle$ such that
$C_\delta\cap S=\emptyset$ for each limit $\delta<\lambda$. As in Abraham
Shelah Solovay \cite[\S 1]{AShS:221} we can have also the squared diamond
sequence.  \QED

\bibliographystyle{literal-plain}
\bibliography{listb,lista,listx}

\shlhetal
\end{document}